\documentclass[a4paper,12pt]{article}
\usepackage[top=2.5cm,bottom=2.5cm,left=2.5cm,right=2.5cm]{geometry}
\usepackage{amsfonts}
\usepackage{latexsym,amsmath,amssymb,amsfonts,epsfig,psfrag,url,graphics,ifpdf,multicol}
\usepackage{color,colordvi,amscd,amsthm} 
\usepackage{tikz}
\usepackage{subfigure}
\usetikzlibrary{decorations.pathreplacing}
\usepackage{verbatim}
\usepackage[numbers,sort&compress]{natbib}
\usepackage{indentfirst,graphics,epsfig}
\usepackage{graphicx}
\usepackage{graphics}
\usepackage{graphicx,psfrag}
\usepackage{graphics}
\usepackage{mathrsfs}
\usepackage{tabularx}
\usepackage{booktabs}
\usepackage{floatrow}
\usepackage{multirow}
\usepackage{graphicx}
\usepackage{makecell}
\usepackage{caption}
\usepackage[ruled]{algorithm2e}
\newfloatcommand{capbtabbox}{table}[][\FBwidth]

\newtheorem{thm}{Theorem}
\newtheorem{lem}{Lemma}

\newtheorem{cl}{Claim}
\newtheorem{rema}{Remark}
\theoremstyle{plain}
\theoremstyle{plain}

\newtheorem{coro}{Corollary}

 \theoremstyle{definition}

\theoremstyle{remark}
\def\pf{\noindent{\bf Proof.\ }}
\def\qed{{\hfill\rule{4pt}{7pt}}}

\allowdisplaybreaks

\numberwithin{subcase}{case}
\numberwithin{subca}{ca}
\numberwithin{subsubcase}{subcase}

\usepackage{indentfirst,subfig}

\usepackage{booktabs} 
\usepackage{pifont}
\usepackage{tikz}

\begin{document}
	\captionsetup[figure]{labelfont={bf},name={Fig.},labelsep=period}

	\begin{center} {\large The minimum number of maximal independent sets in graphs with given order and independence number}
	\end{center}
	\pagestyle{empty}
	
	\begin{center}
		{
			{\small Yuting Tian, Jianhua Tu\footnote{Corresponding author.\\\indent \ \  E-mail: tujh81@163.com (J. Tu)}}\\[2mm]
			
			{\small School of Mathematics and Statistics, Beijing Technology and Business University, \\
				\hspace*{1pt} Beijing, P.R. China 100048}\\[2mm]}
		
	\end{center}
	
	\begin{center}
		\begin{abstract}
			
	\vskip 3mm
	Let $MIS(G)$ be the set of all maximal independent sets in a graph $G$, and let $mis(G)=|MIS(G)|$. In this paper, we show that for any tree $T$ with  $n$ vertices and independence number $\alpha$, 
	\[mis(T)\geq f(n-\alpha),\]
	and for any unicyclic graph $G$ with $n$ vertices and independence number $\alpha$,
	\begin{align*}
		mis(G)\geq 
		\begin{cases}
			2, & \text{if} \ n=4\ \text{and}\ \alpha=2,\\
			3, & \text{if} \; \alpha=n-2 \; \text{and} \; n\neq4, \\
			2f(n-\alpha), & \text{if} \; n\geq 5\; \text{and}\; \lceil \frac{n}{2} \rceil \leq \alpha < n-2,\\
			f(n-\alpha+2)-f(n-\alpha-3), &\text{if} \; n\geq 5, \;\text{and}\ n \; \text{is odd}, \; \text{and} \; \alpha = \lfloor \frac{n}{2} \rfloor,
		\end{cases}
	\end{align*}
	where $f(n)$ represent the $n$th Fibonacci number. Moreover, we also show that the above inequalities are sharp.
		
			
				\end{abstract}
	\end{center}

	\baselineskip=0.24in
	\section{Introduction}\label{sec1}
	
In this paper, we only consider simple, finite, and labeled graphs. Let $G$ be a graph. The order of $G$ is the number of vertices in $G$. Let $v$ be a vertex in $G$. The neighbors of $v$ are the vertices adjacent to $v$, and the set of neighbors of $v$ is denoted by $N_G(v)$. Let $N_G[v]=N_G(v)\cup\{v\}$. The degree of $v$, denoted by $d_G(v)$, is the number of vertices in $N_G(v)$. If $d_G(v)=1$, then $v$ is called a leaf of $G$. A vertex that has a degree of at least 2 and is adjacent to a leaf is called a support vertex. If $G$ does not have two vertices with exactly the same set of neighbors, then $G$ is called twin-free. For $S\subseteq V(G)$, the induced subgraph $G[S]$ consists of the vertices in $S$ and all edges between them in $G$. The graph $G-S$ is obtained by deleting all vertices in $S$ from $G$. Furthermore, $G-\{v\}$ is simply written as $G-v$.
	
An independent set of $G$ is a subset of vertices such that no two vertices in the subset are adjacent. A maximal independent set is an independent set that is not contained in any other independent set, and a maximum independent set is an independent set that contains the maximum possible number of vertices. The independence number $\alpha(G)$ of $G$ is the number of vertices in a maximum independent set. Let $MIS(G)$ denote the set of all maximal independent sets in $G$, and $mis(G)=|MIS(G)|$.

Let $P_n$, $C_n$, and $K_{1,n-1}$ be the path, the cycle, and the star with $n$ vertices, respectively. A unicyclic graph is a connected graph with exactly one cycle.
	
The study of the extremal enumeration problem of maximal independent sets originated from the famous question posed by Erd\H{o}s and Moser in the 1960s:
	
	{\it What is the maximum number of maximal independent sets in a general graph with given order? What are the extremal graphs that achieve this maximum number?}
	
\noindent Since then, this type of problem has attracted widespread attention. Researchers have determined the maximum number of maximal independent sets for various graph classes with given order, as well as the extremal graphs that achieve these maximum numbers. These graph classes include general graphs \cite{Miller1960,Moon1965}, connected graphs \cite{Furedi1987,Griggs1988}, trees \cite{Wilf1986,Sagan1988}, trees with a given number of leaves \cite{Taleskii2022}, triangle-free graphs \cite{Hujtera1993}, unicyclic graphs \cite{Koh2008}, bipartite graphs \cite{Liu1994}, graphs with at most one cycle \cite{Jou1997}, graphs with at most $r$ cycles \cite{Sagan2006}, and graphs without an induced triangle matching of size $t+1$ \cite{Palmer2023}, among others.

Over the past few decades, the primary focus of researchers has been on the maximum number of maximal independent sets in various graph classes with given order, with little attention paid to the minimum number. This is mainly because, in some common classes of graphs, determining the minimum number of maximal independent sets is relatively straightforward. For instance, the minimum number of maximal independent sets in trees or connected graphs of order $n\geq2$ is obviously 2, with the star being an obvious extremal graph. 
	
Therefore, to consider the minimum number of maximal independent sets, more restrictions on the graph classes are needed, or additional graph parameters must be specified beyond the graph's order. Following this line of thought, Talevskii and Malyshev \cite{Taletskii2018,Taleskii2022} considered trees with a given number of leaves and trees without twin-leaves, determining the minimum number of maximal independent sets in these two types of trees with given order. Recently, Cambie and Wagner \cite{Cambie2022} have considered twin-free connected graphs, twin-free bipartite graphs, and twin-free trees, determining the minimum number of maximal independent sets in these graph classes with given order and identifying the corresponding extremal graphs. In this paper, we choose another direction; by specifying more graph parameters, we determine the minimum number of maximal independent sets in trees and unicyclic graphs with given order and independence number.

For any tree $T$ of order $n\geq2$, $\lceil\frac{n}{2}\rceil \leq \alpha(T) \leq n-1$, and the upper bound is reached only for the star. Our first main result is the determination of the minimum number of maximal independent sets in trees with given order and independence number. Let $f(n)$	represent the {\bf $n$th Fibonacci number}, then $f(0)=0$, $f(1)=1$, and for $n\geq 2$, $f(n)=f(n-1)+f(n-2)$.

\begin{thm}\label{thm1}
For any tree $T$ of order $n$ and independence number $\alpha$, 
\[mis(T)\geq f(n-\alpha+2),\]
and this inequality is sharp.
\end{thm}

For any unicyclic graph $G$ of order $n\geq 3$, $\lfloor \frac{n}{2} \rfloor \leq \alpha(G)\leq n-2$. Our second main result is the determination of the minimum number of maximal independent sets in unicyclic graphs with given order and independence number.

\begin{thm}\label{thm2}
For any unicyclic graph $G$ of order $n\geq 3$ and independence number $\alpha$,
\begin{align*}
	mis(G)\geq 
	\begin{cases}
		2, & \text{if} \ n=4\ \text{and}\ \alpha=2,\\
		3, & \text{if} \; \alpha=n-2 \; \text{and} \; n\neq4, \\
		2f(n-\alpha), & \text{if} \; n\geq 5\; \text{and}\; \lceil \frac{n}{2} \rceil \leq \alpha < n-2,\\
		f(n-\alpha+2)-f(n-\alpha-3), &\text{if} \; n\geq 5, \;\text{and}\ n \; \text{is odd}, \; \text{and} \; \alpha = \lfloor \frac{n}{2} \rfloor,
	\end{cases}
\end{align*}
and this inequality is sharp.
\end{thm}

\section{Proof of Theorem \ref{thm1}}\label{sec2}

Let 
\[g(n)=f(n+2).\]
Clearly, $g(0)=f(2)=1$, $g(1)=f(3)=2$ and for $n\geq 2$, $g(n)=g(n-1)+g(n-2)$.

	\begin{lem}\label{lem1}
Let $n_1$ and $n_2$ be two natural numbers. We have

	(1) $g(n_1)\cdot g(n_2)\geq g(n_1+n_2)$, 
	
	(2) if $n_1\geq 1$ and $n_2\geq 1$, then $g(n_1)\cdot g(n_2)+g(n_1-1)\cdot g(n_2-1)=g(n_1+n_2+1)$.
\end{lem}

\pf Since it is clear that for $n\geq 0$, $g(0)\cdot g(n)=g(n)$, we only need to prove that when $n_1\geq 1$ and $n_2\geq 1$, the lemma holds. We proceed by induction on $n_1+n_2$. Firstly, we prove that the lemma holds when $n_1=1$ or $n_2=1$. Given the symmetry of $n_1$ and $n_2$, it suffices to prove the lemma holds when $n_1=1$. We have

(1) $g(1)\cdot g(n_2)=2\cdot f(n_2+2)\geq f(n_2+2)+f(n_2+1)=f(n_2+3)=g(n_2+1)$,

(2) $g(1)\cdot g(n_2)+g(0)\cdot g(n_2-1)=2\cdot g(n_2)+g(n_2-1)=g(n_2)+g(n_2+1)=g(n_2+2).$

\noindent Secondly, we prove that the lemma holds when $n_1=2$ or $n_2=2$. Given the symmetry of $n_1$ and $n_2$, it suffices to prove the lemma holds when $n_1=2$. We have

(1) $g(2)\cdot g(n_2)=3\cdot f(n_2+2)\geq 2\cdot f(n_2+2)+f(n_2+1)=f(n_2+4)=g(n_2+2)$,

(2) $g(2)\cdot g(n_2)+g(1)\cdot g(n_2-1)=3\cdot g(n_2)+2\cdot g(n_2-1)=g(n_2)+2\cdot g(n_2+1)=g(n_2+3).$

Let us assume $n_1\geq 3$, $n_2\geq3$, and that the lemma holds for all pairs of positive integers with a sum less than $n_1+n_2$. Thus,
\begin{align*}
(1)\  g(n_1)\cdot g(n_2)&=[g(n_1-1)+g(n_1-2)]\cdot [g(n_2-1)+g(n_2-2)]\\
	&= g(n_1-1)\cdot g(n_2-1)+g(n_1-1)\cdot g(n_2-2)\\
	&\quad +g(n_1-2)\cdot g(n_2-1)+g(n_1-2)\cdot g(n_2-2) \\
	&\geq g(n_1+n_2-2) + 2\cdot g(n_1+n_2-3)+g(n_1+n_2-4)\\ 
	&= g(n_1+n_2-1)+g(n_1+n_2-2) \\
	&=g(n_1+n_2).
\end{align*}

\begin{align*}
(2)\ &g(n_1)\cdot g(n_2)+g(n_1-1)\cdot g(n_2-1)\\
&=[g(n_1-1)+g(n_1-2)]\cdot [g(n_2-1)+g(n_2-2)]+g(n_1-1)\cdot g(n_2-1)\\
&= g(n_1-1)\cdot g(n_2-1)+g(n_1-2)\cdot g(n_2-2)\\
&\quad +g(n_1-1)\cdot g(n_2-2)+g(n_1-2)\cdot g(n_2-1)+g(n_1-1)\cdot g(n_2-1)\\
&=g(n_1+n_2-1)+g(n_1-1)\cdot g(n_2)+g(n_1-2)\cdot g(n_2-1)\\
&=g(n_1+n_2-1)+g(n_1+n_2)\\
&=g(n_1+n_2+1).
\end{align*} 

The proof of Lemma \ref{lem1} is complete. \qed


For natural numbers $a$, $b$, $c$, and $d$, if $a \geq c$ and $a+b \geq c+d$, we say $(a,b)$ majorizes $(c,d)$, written by $(a,b) \succeq (c,d)$.

\begin{lem}\cite{Cambie2022}\label{lem2}
	For natural numbers $a$, $b$, $c$, $d$, $e$, $f$, $g$, and $h$, if $(a,b) \succeq (c,d)$, $(e,f) \succeq (g,h)$, $c \geq d$, and $g \geq h$, then $(ae,bf) \succeq (cg,dh)$.
\end{lem}






It is agreed that if $G$ is a null graph, then $mis(G)=1$. The following two lemmas are straightforward.

\begin{lem}\label{lem3}
Let $G$ be a graph, and let $y$ be a support vertex in $G$. Let $Q$ be the set of leaves that are adjacent $y$, then 
	\[mis(G)=mis(G-Q-y)+mis(G-N_G[y]).\]
\end{lem}

\begin{lem}\label{lem4}
If $G$ is a disconnect graph with $k$ components $G_1$, $G_2$, $\cdots$, $G_k$, then \[mis(G)=\prod^k_{i=1}mis(G_i).\]
\end{lem}


\vskip 3mm

Now, we are in a position to prove Theorem \ref{thm1}.

\noindent {\bf Proof of Theorem \ref{thm1}.} If $n\leq 2$, it is easy to check that $mis(T)\geq g(n-\alpha)$. Therefore, we can assume that $n\geq 3$. Let $p=n-\alpha$. We prove by induction on $p$ that $mis(T)\geq g(p)$. When $p=1$, $T$ is a star and $mis(T)=2=g(p)$. 

Assume that $p\geq 2$, and that the inequality in the theorem holds for all trees where the number of vertices minus the independence number is less than $p$. 

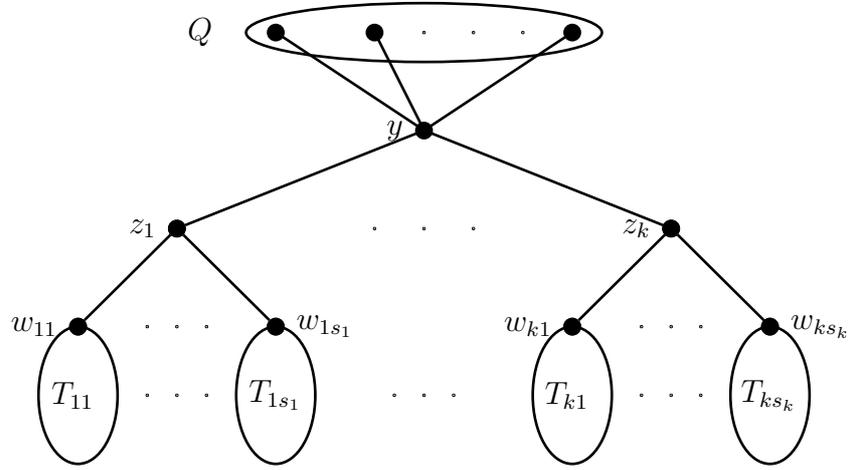
\begin{figure}[H]
	\begin{tikzpicture}
		[line width = 1pt, scale=1.3,
		empty/.style = {circle, draw, fill = white, inner sep=0mm, minimum size=1mm}, full/.style = {circle, draw, fill = black, inner sep=0mm, minimum size=2mm}, full1/.style = {circle, draw, fill = black, inner sep=0mm, minimum size=0.1mm}]

		\draw (4.5,10) ellipse (1.8cm and 0.3cm);
		\node [label=left:$Q$](a1) at (2.6,10) {};
		\node [full1](x1) at (3,10) {};
			\node [full](a1) at (3,10) {};
		\node [full1](x2) at (4,10) {};
		\node [full](a1) at (4,10) {};
		\node [full1](xl) at (6,10) {};
		\node [full](a1) at (6,10) {};
		\node [full1](y) at (4.5,9) {};
		\node [full,label=left:$y$](a1) at (4.5,9) {};
		\node [full1](x3) at (4.5,10) {};
		\node [full1](x4) at (5,10) {};
		\node [full1](x5) at (5.5,10) {};
		\draw (x1) -- (y);
		\draw (x2) -- (y);
		\draw (xl) -- (y);
		
		\node [full](z1) at (2,8) {};
		\node [full,label=left:$z_1$](a1) at (2,8) {};
		\draw (z1) -- (y);
		\node [full](zk) at (7,8) {};
		\node [full,label=left:$z_k$](a1) at (7,8) {};
		\draw (zk) -- (y);
		\node [full1](z3) at (4.5,8) {};
		\node [full1](z4) at (5,8) {};
		\node [full1](z5) at (4,8) {};
		
		\node [full](w11) at (1,7) {};
		\node [full,label=left:$w_{11}$](a1) at (1,7) {};
		\draw (z1) -- (w11);
		\node [full](w1s1) at (3,7) {};
		\node [full,label=right:$w_{1s_1}$](a1) at (3,7) {};
		\draw (z1) -- (w1s1);
		\node [full1](w12) at (1.7,7) {};
		\node [full1](w13) at (2,7) {};
		\node [full1](w14) at (2.3,7) {};
		
		\node [full](wk1) at (6,7) {};
		\node [full,label=left:$w_{k1}$](a1) at (6,7) {};
		\draw (zk) -- (wk1);
		\node [full](wksk) at (8,7) {};
		\node [full,label=right:$w_{ks_k}$](a1) at (8,7) {};
		\draw (zk) -- (wksk);
		\node [full1](wk2) at (6.7,7) {};
		\node [full1](wk3) at (7,7) {};
		\node [full1](wk4) at (7.3,7) {};

		\draw (1,6.3) ellipse (0.4cm and 0.7cm);
		\draw (3,6.3) ellipse (0.4cm and 0.7cm);
		\draw (6,6.3) ellipse (0.4cm and 0.7cm);
		\draw (8,6.3) ellipse (0.4cm and 0.7cm);
		\node [label=left:$T_{11}$](a1) at (1.4,6.3) {};
		\node [label=left:$T_{1s_1}$](a1) at (3.5,6.3) {};
		\node [label=left:$T_{k1}$](a1) at (6.4,6.3) {};
		\node [label=left:$T_{ks_k}$](a1) at (8.5,6.3) {};
		\node [full1](t11) at (1.7,6.3) {};
		\node [full1](t12) at (2.3,6.3) {};
		\node [full1](t13) at (2,6.3) {};
		\node [full1](t21) at (4.5,6.3) {};
		\node [full1](t22) at (4.8,6.3) {};
		\node [full1](t23) at (4.2,6.3) {};
		\node [full1](tk1) at (6.7,6.3) {};
		\node [full1](tk2) at (7,6.3) {};
		\node [full1](tk3) at (7.3,6.3) {};
		
	\end{tikzpicture}
	\caption{The structure of $T$.}
	\label{fig1}
\end{figure}

Let $T$ be a tree of order $n$ and independence number $\alpha$ such that $n-\alpha=p$. Since $p\geq 2$, $T$ is not a star. Let $y$ be a support vertex in $T$, and let $Q$ be the set of leaves that are adjacent to $y$. Let $q=|Q|$. Suppose that
\[N_T(y)=Q\cup\{z_1,z_2,\cdots,z_k\}.\]
For each $i\in\{1,\cdots,k\}$, let $T_i$ be the component in $T-Q-y$ that contains the vertex $z_i$, and let $n_i=|V(T_i)|$, $\alpha_i=\alpha(T_i)$, and $p_i=n_i-\alpha_i$ (clearly, $p_i\geq 1$). Furthermore, let 
\[N_T(z_i)=\{y\}\cup\{w_{i1},\cdots,w_{is_i}\},\]
let $T_{ij}$ be the component in $T_i-z_i$ that contains the vertex $w_{ij}$, and let $n_{ij}=|V(T_{ij})|$, $\alpha_{ij}=\alpha(T_{ij})$, and $p_{ij}=n_{ij}-\alpha_{ij}$. See Figure \ref{fig1}.

By Lemmas \ref{lem3} and \ref{lem4}, 
\begin{align*}
	mis(T)&=mis(T-Q-y)+mis(T-N_T[y])\\&
	=\prod^k_{i=1}mis(T_i)+\prod^k_{i=1}\prod^{s_i}_{j=1}mis(T_{i,j}).
\end{align*}

We show that for every $i\in\{1,\cdots,k\}$, $p_i<p$. Since $|V(T-Q-y)|=n-q-1=n_1+\cdots+n_k$ and $\alpha(T-Q-y)=\alpha-q=\alpha_1+\cdots+\alpha_k$, we have
\[p_1+\cdots+p_k=(n_1-\alpha_1)+\cdots+(n_k-\alpha_k)=(n-q-1)-(\alpha-q)=p-1,\]
which implies that for every $i\in\{1,\cdots,k\}$, $p_i<p$. 

We also have that for any $i$ and any $j$, $p_{ij}\leq p_i$ and $p_i\geq \sum\limits_{j=1}^{s_i}p_{ij}\geq p_i-1$. This is because $\alpha_i-1\leq\alpha(T_i-z_i)\leq \alpha_i$, and
\[p_i=(n_i-1)-(\alpha_i-1)\geq \sum\limits_{j=1}^{s_i}p_{ij}=|V(T_i-z_i)|-\alpha(T_i-z_i)\geq (n_i-1)-\alpha_i=p_i-1.\]

By the inductive hypothesis, for any $i$ and any $j$, $mis(T_i)\geq g(p_i)$ and $mis(T_{ij})\geq g(p_{ij}).$ By Lemma \ref{lem1}, for each $i$, 
\[\prod_{j=1}^{s_i}mis(T_{ij})\geq \prod_{j=1}^{s_i}g(p_{ij}) \geq g\left(\sum^{s_i}_{j=1}p_{ij}\right) \geq g(p_i-1), \]
which implies
\[\left(mis(T_i),\ \prod_{j=1}^{s_i}mis(T_{ij})\right) \succeq \left(g(p_i),\ g(p_i-1)\right).\]

By Lemma \ref{lem2}, 
\[\left(\prod_{i=1}^2mis(T_i),\ \prod_{i=1}^2\prod_{j=1}^{s_i}mis(T_{ij})\right) \succeq \left(g(p_1)g(p_2),\ g(p_1-1)g(p_2-1)\right).\]
which implies that 
\begin{align*}
&\prod_{i=1}^2mis(T_i)+\prod_{i=1}^2\prod_{j=1}^{s_i}mis(T_{ij})\\
&\geq g(p_1)g(p_2)+g(p_1-1)g(p_2-1)\\
&=g(p_1+p_2+1)\\
&=g(p_1+p_2)+g(p_1+p_2-1).
\end{align*}

On the other hand, $\prod_{i=1}^2mis(T_i)\geq g(p_1)\cdot g(p_2)=g(p_1+p_2)$. Thus,
\[\left(\prod_{i=1}^2mis(T_i),\ \prod_{i=1}^2\prod_{j=1}^{s_i}mis(T_{ij})\right) \succeq \left(g(p_1+p_2),\ g(p_1+p_2-1)\right).\]

Furthermore, by Lemma \ref{lem2}, 
\[\left(\prod_{i=1}^3mis(T_i), \ \prod_{i=1}^3\prod_{j=1}^{s_i}mis(T_{ij})\right) \succeq \left(g(p_1+p_2)g(p_3), \ g(p_1+p_2-1)g(p_3-1)\right),\]
which implies that 
\begin{align*}
	&\prod_{i=1}^3mis(T_i)+\prod_{i=1}^3\prod_{j=1}^{s_i}mis(T_{ij})\\
	&\geq g(p_1+p_2)g(p_3)+g(p_1+p_2-1)g(p_3-1)\\
	&=g(p_1+p_2+p_3+1)\\
	&=g(p_1+p_2+p_3)+g(p_1+p_2+p_3-1).
\end{align*}

By repeating this process, we can finally obtain that
\[\prod^{k}_{i=1}mis(T_i)+\prod^{k}_{i=1}\prod^{s_i}_{j=1}mis(T_{ij}) \geq g\left(\sum^k_{i=1}p_i\right)+g\left(\sum^k_{i=1}p_i-1\right).\]	
Since $\sum^k_{i=1}p_i=p-1$, 
\begin{align*}
	mis(T)&=\prod^{k}_{i=1}mis(T_i)+\prod^{k}_{i=1}\prod^{s_i}_{j=1}mis(T_{ij})\\
	&\geq g\left(\sum^k_{i=1}p_i\right)+g\left(\sum^k_{i=1}p_i-1\right)\\
	&=g(p-1)+g(p-2)\\
	&=g(p).
\end{align*}

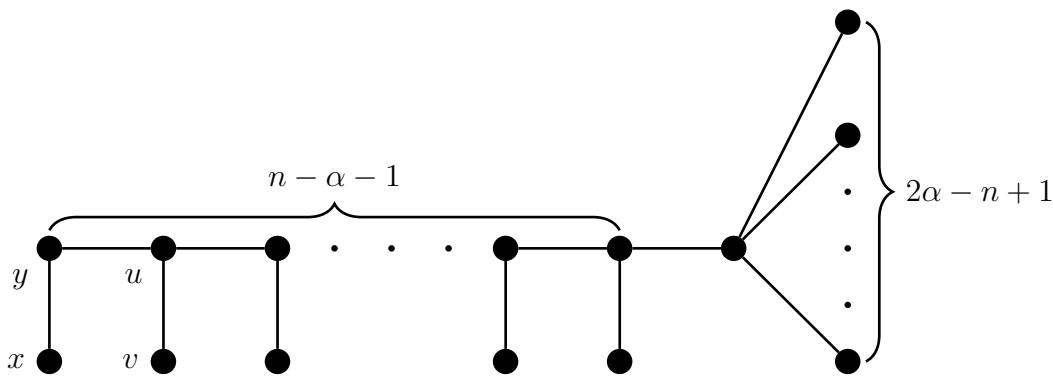
\begin{figure}[H]
	  	\begin{tikzpicture}
	  		[line width = 1pt, scale=1.5,
	  		empty/.style = {circle, draw, fill = white, inner sep=0mm, minimum size=2mm}, full/.style = {circle, draw, fill = black, inner sep=0mm, minimum size=3mm}, full1/.style = {circle, draw, fill = black, inner sep=0mm, minimum size=0.5mm}]
	  		\node [full,label=below left:$y$](a) at (1,1) {};
	  		\node [full,,label=left:$x$](a1) at (1,0) {};
	  		\node [full,label=below left:$u$](b) at (2,1) {};
	  		\node [full,label=left:$v$](b1) at (2,0) {};
	  		\node [full](c) at (3,1) {};
	  		\node [full](c1) at (3,0) {};
	  		\node [full](d) at (5,1) {};
	  		\node [full](d1) at (5,0) {};
	  		\node [full](e) at (6,1) {};
	  		\node [full](e1) at (6,0) {};
	  		\node [full](f) at (7,1) {};
	  		\node [full](f1) at (8,3) {};
	  		\node [full](f2) at (8,2) {};
	  		\node [full1](f3) at (8,1.5) {};
	  		\node [full1](f4) at (8,1) {};
	  		\node [full1](f5) at (8,0.5) {};
	  		\node [full](f6) at (8,0) {};
	  		\node [full1](u) at (3.5,1) {};
	  		\node [full1](v) at (4,1) {};
	  		\node [full1](w) at (4.5,1) {};
	  		
	  		\draw (a) -- (b);
	  		\draw (b) -- (c);
	  		\draw (d) -- (e);
	  		\draw (e) -- (f);
	  		\draw (a) -- (a1);
	  		\draw (b) -- (b1);
	  		\draw (c) -- (c1);
	  		\draw (d) -- (d1);
	  		\draw (e) -- (e1);
	  		\draw (f) -- (f1);
	  		\draw (f) -- (f2);
	  		\draw (f) -- (f6);
	  		\draw[decorate,decoration={brace,amplitude=10pt},xshift=-4pt,yshift=0pt]
	  		(8.3,3) -- (8.3,0) node [black,midway,xshift=1.5cm] {$2\alpha-n+1$};
	  		\draw[decorate,decoration={brace,amplitude=10pt},xshift=0pt,yshift=-4pt]
	  		(1,1.3) -- (6,1.3) node [black,midway,yshift=0.7cm] {$n-\alpha-1$};
	  	\end{tikzpicture}
	  	\caption{A tree $T$ of order $n$ and independence number $\alpha$ satisfying $mis(T)=g(n-\alpha)$.}
	  	\label{fig2}
	  \end{figure}

Next, we prove that the inequality in the theorem is sharp. For any integers $n$ and $\alpha$ such that $n\geq 2$ and $\lceil\frac{n}{2}\rceil \leq \alpha \leq n-1$, we can construct a tree $T_{n,\alpha}$ of order $n$ and independence number $\alpha$ that is shown in Figure \ref{fig2}.

We prove by induction on $n-\alpha$ that the tree $T_{n,\alpha}$ has $g(n-\alpha)$ maximal independent sets. 
When $n-\alpha=1$, the tree $T_{n,\alpha}$ is a star and $mis(T)=2=g(1).$ 

Assume that $k\geq 2$, and that when $n-\alpha<k$, the tree $T_{n-\alpha}$ has $g(n-\alpha)$ maximal independent sets. When $n-\alpha=k$,
\begin{align*}
mis(T_{n,\alpha})&=mis(T_{n,\alpha}-\{x,y\})+mis(T_{n,\alpha}-\{x,y,u,v\})\\
&=mis(T_{n-2,\alpha-1})+mis(T_{n-4,\alpha-2})\\
&=g(n-\alpha-1)+g(n-\alpha-2)\\
&=g(n-\alpha),
\end{align*}
where $x$, $y$, $u$, and $v$ are four vertices in $T$, shown in Figure \ref{fig2}.

The proof of Theorem \ref{thm1} is complete. \qed

\begin{coro}\label{coro1}
For any forest $F$ of order $n$ and independence number $\alpha$,
\[mis(F)\geq g(n-\alpha),\]
and this inequality is sharp.
\end{coro}

\pf Let $F_1$, $F_2$, $\cdots$, $F_k$ be components of $F$. For each $i\in \{1,\cdots,k\}$, let $n_i=|V(F_i)|$ and $\alpha_i=\alpha(F_i)$. Thus,
\[n_1+\cdots+n_k=n,\ \alpha_1+\cdots+\alpha_k=\alpha.\]

By Lemmas \ref{lem1} and \ref{lem4}, and Theorem \ref{thm1},
\[mis(F)\geq \prod\limits_{i=1}^{k}g(n_i-\alpha_i)\geq g(n-\alpha).\]
It can also be seen from the graph shown in Figure \ref{fig2} that for forests, this inequality is sharp. \qed

	  
\section{Proof of Theorem \ref{thm2}}\label{sec3}

For any integer $n\geq0$, let 
\[h(n)=2f(n).\] 
We have that $h(0)=0$, $h(1)=2$, and for $n\geq 2$, $h(n)=h(n-1)+h(n-2)$.

For any integer $n\geq 3$, let 
\[\ell(n)=f(n+2)-f(n-3).\]
We have that $\ell(3)=5$, $\ell(4)=7$, and for $n\geq 5$, 
\begin{align*}
	\ell(n)&=f(n+2)-f(n-3)\\
	&=f(n+1)+f(n)-f(n-4)-f(n-5)\\
	&=\ell(n-1)+\ell(n-2).
\end{align*}
Moreover, \[\ell(n)=2f(n)+f(n-2)=h(n)+f(n-2)\geq h(n).\]

\begin{lem}\label{lem5}
For any integers $n_1\geq 2$ and $n_2\geq 0$,

\[h(n_1)\cdot g(n_2)\geq h(n_1+n_2).\]
\end{lem}

\pf Prove by induction on $n_1+n_2$.

When $n_1=2$, $h(2)\cdot g(n_2)=2g(n_2)=2f(n_2+2)=h(n_2+2).$

When $n_1=3$, $h(3)\cdot g(n_2)=4g(n_2)=4f(n_2+2)\geq 2f(n_2+3)=h(n_2+3).$

When $n_2=0$, $h(n_1)\cdot g(0)=h(n_1)=h(n_1+0).$

When $n_2=1$, $h(n_1)\cdot g(1)=2h(n_1)\geq h(n_1+1).$

Let us assume $n_1\geq 4$, $n_2\geq2$, and that the lemma holds for all pairs of positive integers with a sum less than $n_1+n_2$. Thus,
\begin{align*}
	h(n_1)\cdot g(n_2)&=[h(n_1-1)+h(n_1-2)]\cdot [g(n_2-1)+g(n_2-2)]\\
	&= h(n_1-1)\cdot g(n_2-1)+h(n_1-1)\cdot g(n_2-2)\\
	&\quad +h(n_1-2)\cdot g(n_2-1)+h(n_1-2)\cdot g(n_2-2) \\
	&\geq h(n_1+n_2-2) + 2\cdot h(n_1+n_2-3)+h(n_1+n_2-4)\\ 
	&= h(n_1+n_2-1)+h(n_1+n_2-2) \\
	&=h(n_1+n_2).
\end{align*}

The proof of Lemma \ref{lem5} is complete. \qed


\begin{lem}\label{lem6}
	For any integers $n_1\geq 3$ and $n_2\geq 0$,
	
	\[\ell(n_1)\cdot g(n_2)\geq \ell(n_1+n_2).\]
\end{lem}

\pf Since $\ell(3)=5$, $\ell(4)=7$ and for $n\geq 5$, $\ell(n)=\ell(n-1)+\ell(n-2)$, similar to the proof of Lemma \ref{lem5}, we can prove that Lemma \ref{lem6} holds. \qed


\begin{lem}\label{lem7}(\cite{Furedi1987})
$mis(C_3)=3$, $mis(C_4)=2$, $mis(C_5)=5$, and for $n\geq 6$, \[mis(C_n)=mis(C_{n-2})+mis(C_{n-3}).\] 
\end{lem}






\begin{lem}\label{lem9}
	If $n\geq 5$, then 
	\[mis(C_n)\geq \ell(\lfloor \frac{n+1}{2} \rfloor).\]
\end{lem}

\pf Prove by induction on $n$.

When $n=5$, $mis(C_5)=5=\ell(3)$. When $n=6$, $mis(C_6)=5=\ell(3)$. When $n=7$, $mis(C_7)=7=\ell(4).$

Assume that $n\geq 8$, and for any integer $5\leq k\leq n-1$, $mis(C_k)\geq \ell(\lfloor \frac{k+1}{2} \rfloor)$. 
Now, by Lemma \ref{lem7} and inductive hypothesis, if $n$ is even
\begin{align*}
	mis(C_n)&=mis(C_{n-2})+mis(C_{n-3})\\
	&\geq \ell(\lfloor \frac{n-1}{2} \rfloor)+\ell(\lfloor \frac{n-2}{2} \rfloor)\\
	&=\ell(\frac{n}{2}-1)+\ell(\frac{n}{2}-1)\\
	&\geq\ell(\frac{n}{2})\\
	&=\ell(\lfloor \frac{n+1}{2} \rfloor),
\end{align*}
if $n$ is odd, 
\begin{align*}
	mis(C_n)&=mis(C_{n-2})+mis(C_{n-3})\\
	&\geq \ell(\lfloor \frac{n-1}{2} \rfloor)+\ell(\lfloor \frac{n-2}{2} \rfloor)\\
	&=\ell(\frac{n-1}{2})+\ell(\frac{n-3}{2})\\
	&=\ell(\frac{n+1}{2})\\
	&=\ell(\lfloor \frac{n+1}{2} \rfloor).
\end{align*}
The proof of Lemma \ref{lem9} is complete. \qed


\begin{lem}\label{lem10}
	Let $n_1\geq 3$ and $n_2\geq 0$ be two integers, then 
	\[\ell(n_1)\cdot 2^{n_2} \geq \ell(n_1+n_2).\]
\end{lem}

\pf Prove by induction on $n_1+n_2$.

When $n_1=3$, $\ell(3)\cdot 2^{n_2}=5\cdot 2^{n_2}\geq\ell(n_2+3)$. 

When $n_1=4$, $\ell(4)\cdot 2^{n_2}=7\cdot 2^{n_2}\geq\ell(n_2+4)$. 

When $n_2=0$, $\ell(n_1)=\ell(n_1+0)$.

Let us assume $n_1\geq 5$, $n_2\geq1$, and that the lemma holds for all pairs of positive integers with a sum less than $n_1+n_2$. Thus,
\begin{align*}
    \ell(n_1)\cdot 2^{n_2}&=\left(\ell(n_1-1)+\ell(n_1-2)\right)\cdot 2^{n_2}\\
	&=\ell(n_1-1)\cdot 2^{n_2}+\ell(n_1-2)\cdot 2^{n_2}\\
	&\geq\ell(n_1+n_2-1)+\ell(n_1+n_2-2)\\
	&=\ell(n_1+n_2).
\end{align*}

The proof of Lemma \ref{lem10} is complete. \qed


\vskip 3mm

Now, we are in a position to prove Theorem \ref{thm2}.

\noindent {\bf Proof of Theorem \ref{thm2}.} Let $G$ be a unicyclic graph of order $n$ and independence number $\alpha$.
Let $C$ be the cycle in $G$. For two vertices $u,v\in V(G)$, we define the distance between $u$ and $v$, denoted by $d(u,v)$, as the length of a shortest path from $u$ to $v$ in $G$. For a nonempty set of vertices $S\subset V(G)$, define $d(u,S)=\min\{d(u,v):v\in S\}$ as the distance from $u$ to $S$. Let 
\[d=\max\{d(u,V(C)):u\in V(G)\}.\]
If $d>0$, then the vertices that are farthest from $C$ must all be leaves. We denote one of these leaves as $x$, and the adjacent vertex of $x$ as $y$.
Furthermore, let $Q$ be the set of leaves in $G$ adjacent to $y$, and let $q=|Q|$. 

We prove by induction on $n-\alpha$ that the following inequality holds.
\begin{align}\label{inequ1}
	mis(G)\geq t(n,\alpha)=
	\begin{cases}
		2 & \text{if} \ n=4\ \text{and}\ \alpha=2,\\
		3 & \text{if} \; \alpha=n-2 \; \text{and} \; n\neq4, \\
		h(n-\alpha) & \text{if} \; n\geq 5\; \text{and}\; \lceil \frac{n}{2} \rceil \leq \alpha < n-2,\\
		\ell(n-\alpha) &\text{if} \; n\geq 5, \;\text{and}\ n \; \text{is odd}, \; \text{and} \; \alpha = \lfloor \frac{n}{2} \rfloor.
	\end{cases}
\end{align}

By direct calculations, for any unicyclic graph $G$ of order $n\leq 5$ and independence number $\alpha$, $mis(G)\geq t(n,\alpha)$. Thus, we suppose that $n\geq6$.

When $n-\alpha=2$, either $C$ is of order $3$, with one or two vertices in $C$ connected to a total of $n-3$ leaves, or $C$ is of order $4$, with one vertex or two non-adjacent vertices in $C$ connected to a total of $n-4$ leaves. A few simple calculations show that $mis(G)\geq 3.$ 

\begin{rema}\label{rem1}
When $n-\alpha=2$,
\[t(n,\alpha)\geq 2=h(2).\]
\end{rema}

Assume that $n-\alpha\geq 3$, and that Inequality (\ref{inequ1}) holds for all unicyclic graphs where the number of vertices minus the independence number is less than $n-\alpha$.

We distinguish the following two cases.

\noindent {\bf Case 1.} $\lceil \frac{n}{2} \rceil \leq \alpha < n-2 $.

\noindent{\bf Subcase 1.1.} $d=0$.

In this subcase, $G$ is a cycle of order $n$ and $n$ is even. Thus, $\alpha=n/2$ and by Lemma \ref{lem9},
\[mis(G)\geq \ell(\lfloor\frac{n+1}{2}\rfloor)\geq h(\lfloor\frac{n+1}{2}\rfloor)=h(n/2)=h(n-\alpha).\]

\noindent{\bf Subcase 1.2.} $d=1$.

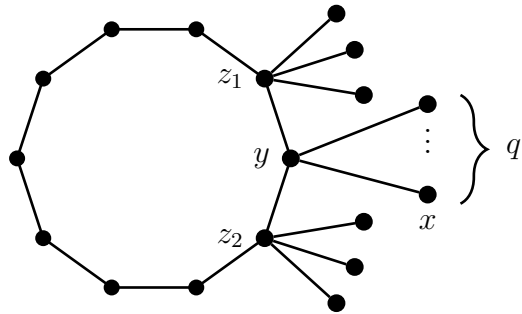
\begin{figure}[H]
		\begin{tikzpicture}
			[line width = 1pt, scale=1.2,
			empty/.style = {circle, draw, fill = white, inner sep=0mm, minimum size=1mm}, full/.style = {circle, draw, fill = black, inner sep=0mm, minimum size=2mm}, full1/.style = {circle, draw, fill = black, inner sep=0mm, minimum size=2mm}]
			
			\foreach \x in {1,...,10} {
				\draw ({360/10 * (\x - 1)}:1.5) -- ({360/10 * \x}:1.5);
				\fill ({360/10 * \x}:1.5) circle (2.5pt);
			}
			\node [full1](x1) at (3,0.6) {};
			\node [full1,label=below:$x$](x2) at (3, -0.4) {};
			\node [full1,label=left:$y$] (y) at ({360/10 * (10)}:1.5) {};
			\node [label=above:$\vdots$] at (3,-0.2) {};
			\draw (y)--(x1);
			\draw (y)--(x2);
			
			\node [full1](z) at (2.2,1.2) {};
			\node [full1](z1) at (2,1.6) {};
			\node [full1](z2) at (2.3,0.7) {};
			\node [full1,label=left:$z_1$](a) at ({360/10 * (1)}:1.5) {};
			\draw (a)--(z);
			\draw (a)--(z1);
			\draw (a)--(z2);
			
			\node [full1,label=left:$z_2$](b) at ({360/10 * (9)}:1.5) {};
			\node [full1](w) at (2.2,-1.2) {};
			\node [full1](w1) at (2,-1.6) {};
			\node [full1](w2) at (2.3,-0.7) {};
			\draw (b)--(w);
			\draw (b)--(w1);
			\draw (b)--(w2);
			
			\draw[decorate,decoration={brace,amplitude=10pt},xshift=-4pt,yshift=0pt]
			(3.5,0.7) -- (3.5,-0.5) node [black,midway,xshift=0.7cm] {$q$};
		\end{tikzpicture}
		\caption{Subcase 1.2. $d=1$}
		\label{fig3}
	\end{figure}

In this subcase, $y\in V(C)$. See Figure \ref{fig3}. We have
\[mis(G)=mis(G-y-Q)+mis(G-N_G[y]).\]
Moreover, $(G-y-Q)$ is a tree of order $n-(q+1)$ and independence number $\alpha-q$, and $G-N_G[y]$ is a forest of order $n-(q+3)$ and independence number at most $\alpha-q$.
By Theorem \ref{thm1} and Corollary \ref{coro1}, 
	\begin{align*}
		mis(G) &=mis(G-y-Q)+mis(G-N_G[y]) \\
		&\geq g\left(n-(q+1)-(\alpha-q)\right)+g\left(n-(q+3)-(\alpha-q)\right) \\
		&\geq g(n-\alpha-1)+g(n-\alpha-3) \\
		&=f(n-\alpha+1)+f(n-\alpha-1)\\ 
		&\geq 2f(n-\alpha)=h(n-\alpha).
\end{align*}
 
\noindent{\bf Subcase 1.3.} $d=2$.

\begin{figure}[H]
		\begin{tikzpicture}
			[line width = 1pt, scale=1.2,
			empty/.style = {circle, draw, fill = white, inner sep=0mm, minimum size=2mm}, full/.style = {circle, draw, fill = black, inner sep=0mm, minimum size=2mm}, full1/.style = {circle, draw, fill = black, inner sep=0mm, minimum size=2mm}]
			
			\foreach \x in {1,...,10} {
				\draw ({360/10 * (\x - 1)}:1.5) -- ({360/10 * \x}:1.5);
				\fill ({360/10 * \x}:1.5) circle (2.5pt);
			}
			
			\node [full1](y1) at (3,0.6) {};
			\node [full1](y2) at (3,-0.6) {};
			\node [full1,label=right:$y$] (y) at (3,0) {};
			\node [full1,label=left:$z$] (z) at ({360/10 * (10)}:1.5) {};
			\draw (z)--(y);
			\draw (z)--(y1);
			\draw (z)--(y2);
			\node [full1](x1) at (4,0.6) {};
			\node [full1,label=below:$x$](x2) at (4,-0.6) {};
			\node [label=above:$\vdots$] at (4,-0.4) {};
			\draw (y)--(x1);
			\draw (y)--(x2);
			
			\node [full1](u) at (2.2,1.2) {};
			\node [full1](u1) at (2,1.6) {};
			\node [full1](u12) at (3,1.6) {};
			\node [full1](u2) at (2.3,0.7) {};
			\node [full1](a) at ({360/10 * (1)}:1.5) {};
			\draw (a)--(u);
			\draw (a)--(u1);
			\draw (a)--(u2);
			\draw (u1)--(u12);
			
			\node [full1](b) at ({360/10 * (9)}:1.5) {};
			\node [full1](w) at (2.2,-1.2) {};
			\node [full1](w1) at (2.3,-0.7) {};
			\node [full1](w2) at (2,-1.6) {};
			\draw (b)--(w);
			\draw (b)--(w1);
			\draw (b)--(w2);
			\node [full1](w02) at (3,-1.2) {};
			\node [full1](w22) at (3,-1.6) {};
			\draw (w2)--(w22);
			\draw (w)--(w02);
			
			\draw[decorate,decoration={brace,amplitude=10pt},xshift=-4pt,yshift=0pt]
			(4.5,0.7) -- (4.5,-0.7) node [black,midway,xshift=0.7cm] {$q$};
		\end{tikzpicture}
		\caption{Subcase 1.3. $d=2$}
		\label{fig4}
	\end{figure}
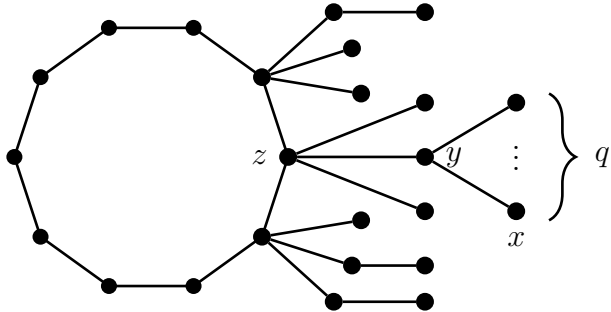

In this subcase, there is a vertex, say $z$, in $C$ that is adjacent to $y$ and $d_G(y)=q+1$. See Figure \ref{fig4}. Moreover, $(G-y-Q)$ is a unicyclic graph of order $n-(q+1)$ and independence number $\alpha-q$, and $G-N_G[y]$ is a forest of order $n-(q+2)$ and independence number at most $\alpha-q$.

If $n-(q+1)$ is odd, and $\alpha-q=\lfloor\frac{n-(q+1)}{2}\rfloor$, and $n-\alpha=3$, then 
$q=n-4$ and $G$ is a graph obtained by identifying a vertex of a triangle with a leaf of the star $K_{1,n-3}$. Thus, 
\[mis(G)=5>4=h(3).\]

If $n-(q+1)$ is odd, and $\alpha-q=\lfloor\frac{n-(q+1)}{2}\rfloor$, and $n-\alpha\geq4$, then by inductive hypothesis, 
	\[mis(G-y-Q)\geq \ell(n-(q+1)-(\alpha-q))=\ell(n-\alpha-1)\geq h(n-\alpha-1).\]
If $n-(q+1)$ is even or $\alpha-q>\lfloor\frac{n-(q+1)}{2}\rfloor$, then by inductive hypothesis, 
\[mis(G-y-Q)\geq h(n-(q+1)-(\alpha-q))=h(n-\alpha-1).\]
In conclusion,
\begin{align*}
	mis(G) &=mis(G-y-Q)+mis(G-N_G[y]) \\
	&\geq h(n-\alpha-1)+g\left(n-(q+2)-(\alpha-q)\right) \\
	&\geq h(n-\alpha-1)+g(n-\alpha-2)\\
	&=2f(n-\alpha-1)+f(n-\alpha)\\
	&\geq 2f(n-\alpha)=h(n-\alpha).
\end{align*}

\noindent{\bf Subcase 1.4.} $d\geq 3$.

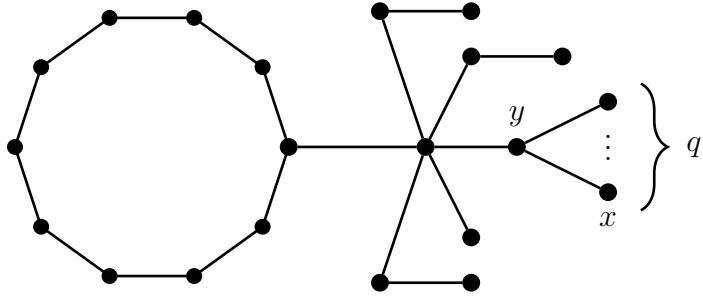
\begin{figure}[H]
		\begin{tikzpicture}
			[line width = 1pt, scale=1.2,
			empty/.style = {circle, draw, fill = white, inner sep=0mm, minimum size=2mm}, full/.style = {circle, draw, fill = black, inner sep=0mm, minimum size=2mm}, full1/.style = {circle, draw, fill = black, inner sep=0mm, minimum size=2mm}]
			
			\foreach \x in {1,...,10} {
				\draw ({360/10 * (\x - 1)}:1.5) -- ({360/10 * \x}:1.5);
				\fill ({360/10 * \x}:1.5) circle (2.5pt);
			}
			\node [full1](u) at ({360/10 * 10}:1.5) {};
			\node [full1](z1) at (2.5,1.5) {};
			\node [full1](z12) at (3.5,1.5) {};
			\node [full1](z2) at (3.5,1) {};
			\node [full1](z22) at (4.5,1) {};
			\node [full1](z3) at (2.5,-1.5) {};
			\node [full1](z32) at (3.5,-1.5) {};
			\node [full1](z4) at (3.5,-1) {};
			\node [full1](z) at (3,0) {};
			\draw (u)--(z);
			\draw (z)--(z1);
			\draw (z)--(z2);
			\draw (z)--(z3);
			\draw (z)--(z4);
			\draw (z1)--(z12);
			\draw (z2)--(z22);
			\draw (z3)--(z32);

			\node [full1,label=above:$y$] (y) at (4,0) {};
			
				\node [label=above:$\vdots$] at (5,-0.4) {};
			\node [full1](x1) at (5,0.5) {};
			\node [full1,label=below:$x$](x2) at (5,-0.5) {};
			\draw (z)--(y);
			\draw (y)--(x1);
			\draw (y)--(x2);
			
				\draw[decorate,decoration={brace,amplitude=10pt},xshift=-4pt,yshift=0pt]
			(5.5,0.7) -- (5.5,-0.7) node [black,midway,xshift=0.7cm] {$q$};			
		\end{tikzpicture}
		\caption{Subcase 1.4. $d \geq 3$}
		\label{fig5}
	\end{figure}

In this subcase, $V(C)\cap N_G(y)=\emptyset$ and $d_G(y)=q+1$. See Figure \ref{fig5}. Moreover, $(G-y-Q)$ is a unicyclic graph of order $n-(q+1)$ and independence number $\alpha-q$, and $G-N_G[y]$ is a disjoint union of a unicyclic graph $G'$ and a forest $F$. Let $n_1=|V(G')|$, $\alpha_1=\alpha(G')$, $n_2=|V(F)|$, and $\alpha_2=\alpha(F)$. Thus,
\[n_1+n_2=n-(q+2),\ \alpha_1+\alpha_2\leq \alpha-q.\]

By Remark \ref{rem1}, and Lemma \ref{lem5}, and the fact that for any integer $m\geq 3$, $\ell(m)\geq h(m)$, we have 
\[mis(G-N_G[y])\geq h(n_1-\alpha_1)\cdot g(n_2-\alpha_2)\geq h(n_1+n_2-(\alpha_1+\alpha_2))\geq h(n-\alpha-2),\]
and
\begin{align*}
	mis(G) &=mis(G-y-Q)+mis(G-N_G[y]) \\
	&\geq h\left(n-(q+1)-(\alpha-q)\right)+h(n-\alpha-2) \\
	&\geq h(n-\alpha-1)+h(n-\alpha-2)\\
	&=h(n-\alpha).
\end{align*}

\vskip 3mm

\noindent {\bf Case 2.} $n\geq 7$, and $n$ is odd, and $\alpha=\lfloor\frac{n}{2} \rfloor$.

We need to prove that
\[mis(G)\geq \ell(n-\alpha)=\ell(\frac{n+1}{2}).\]

Firstly, we prove that there is no vertex in $G$ that is adjacent to at least two leaves. Suppose, to the contrary, that there exists a vertex, say $v$, in $G$ which is adjacent to $k$ leaves $v_1$, $v_2$, $\cdots$, $v_k$, where $k\geq 2$.
Thus,
\begin{align*}
\alpha(G)&=\alpha(G-\{v,v_1,\cdots,v_k\})+k\\
&\geq \lfloor\frac{n-(k+1)}{2}\rfloor+k\\
&= \lfloor\frac{n-1+k}{2}\rfloor\\
&>\lfloor\frac{n}{2}\rfloor,
\end{align*}
which is a contradiction.

Secondly, we prove that $d\neq 1$. Suppose, to the contrary, that $d=1$. Then $G-\{x,y\}$ is a tree of order $n-2$, which implies that $\alpha(G-\{x,y\})\geq \lceil\frac{n-2}{2}\rceil$. Thus, 
\[\alpha(G)=\alpha(G-\{x,y\})+1\geq \lceil\frac{n-2}{2}\rceil+1=\lceil\frac{n}{2}\rceil>\lfloor\frac{n}{2}\rfloor,\]
which is a contradiction.

Thirdly, we show that the order of $C$ is odd by proving the following claim.

\begin{cl}\label{cl1}
For a unicyclic graph $H$ with $m$ vertices, if the order of the cycle in $H$ is even, then $\alpha(H)\geq \lceil\frac{m}{2}\rceil$.
\end{cl}

\noindent{\sl Proof of Claim \ref{cl1}}. Prove by induction on $m$.

When $m=4$, $H$ is a cycle and $\alpha(H)=2=4/2.$

Assume that $m\geq 5$, and for any integer $k\leq m-1$, a unicyclic graph of order $k$ that contains an even cycle has independence number at least $\lceil\frac{k}{2}\rceil$.

Let $C'$ be the cycle in $H$. Let 
\[d'=\max\{d(u,V(C')):u\in V(H)\}.\]
If $d'>0$, we denote one of leaves that are farthest from $C'$ as $x'$, and the adjacent vertex of $x'$ as $y'$.
Furthermore, let $Q'$ be the set of leaves in $H$ adjacent to $y'$, and let $q'=|Q'|$.

If $d'=0$, then $H$ is a cycle of order $m$ and has independence number $m/2$.

If $d'=1$, then $H-y'-Q'$ is a tree, which implies that $\alpha(H-y'-Q')\geq \lceil\frac{m-(q'+1)}{2}\rceil$. Thus, 
\[\alpha(H)=\alpha(H-y'-Q')+q'\geq \lceil\frac{m-(q'+1)}{2}\rceil+q'\geq \lceil\frac{m}{2}\rceil.\]

If $d'\geq 2$, then $H-y'-Q'$ is a unicyclic graph that contains an even cycle. By inductive hypothesis, 
\[\alpha(H)=\alpha(H-y'-Q')+q'\ge\lceil \frac{m-(q'+1)}{2}\rceil+q'\geq \lceil\frac{m}{2}\rceil.\] 

The proof of Claim \ref{cl1} is complete. \qed\\

We distinguish the following three subcases.

\noindent {\bf Subcase 2.1.} $d=0$.

In this subcase, $G$ is an odd cycle. By Lemma \ref{lem9},
\[mis(G)=mis(C_n)\geq \ell(\frac{n+1}{2}).\]

\noindent {\bf Subcase 2.2.} $d=2$.
	
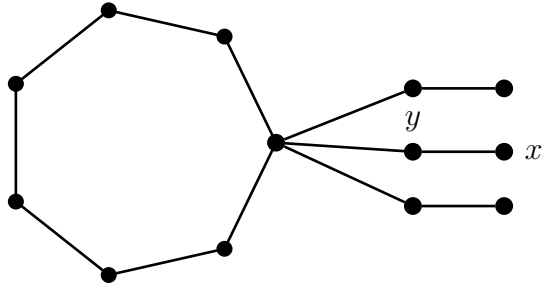
\begin{figure}[H]
		\begin{tikzpicture}
			[line width = 1pt, scale=1.2,
			empty/.style = {circle, draw, fill = white, inner sep=0mm, minimum size=2mm}, full/.style = {circle, draw, fill = black, inner sep=0mm, minimum size=2mm}, full1/.style = {circle, draw, fill = black, inner sep=0mm, minimum size=2mm}]
			
			\foreach \x in {1,...,7} {
				\draw ({360/7 * (\x - 1)}:1.5) -- ({360/7 * \x}:1.5);
				\fill ({360/7 * \x}:1.5) circle (2.5pt);
			}
			
			\node [full1](y1) at (3,0.6) {};
			\node [full1](y2) at (3,-0.7) {};
			\node [full1,label=above:$y$] (y) at (3,-0.1) {};
			\node [full1] (z) at ({360/7 * (7)}:1.5) {};
			\draw (z)--(y);
			\draw (z)--(y1);
			\draw (z)--(y2);
			\node [full1,label=right:$x$] (x) at (4,-0.1) {};
			\node [full1](y12) at (4,0.6) {};
			\node [full1](y22) at (4,-0.7) {};
			\draw (y)--(x);
			\draw (y1)--(y12);
			\draw (y2)--(y22);
						
		\end{tikzpicture}
		\caption{Subcase 2.2. $d=2$}
		\label{fig6}
	\end{figure}

In this subcase, there is a vertex in $C$ that is adjacent to $y$ and $d_G(y)=2$. See Figure \ref{fig6}. We have
\[mis(G)=mis(G-\{x,y\})+mis(G-N_G[y]).\]
Moreover, $G-\{x,y\}$ is a unicyclic graph of order $n-2$ and independence number $\lfloor \frac{n}{2}\rfloor-1$, and $G-N_G[y]$ is a forest of order $n-3$ and independence number at most $\lfloor \frac{n}{2}\rfloor-1$. By inductive hypothesis, 
\begin{align*}
	mis(G)&=mis(G-\{x,y\})+mis(G-N_G[y])\\
	&\geq \ell(\frac{n-1}{2})+g(n-3-\frac{n-3}{2})\\
	&\geq\ell(\frac{n+1}{2}).
\end{align*}

\noindent {\bf Subcase 2.3.} $d\geq 3$.

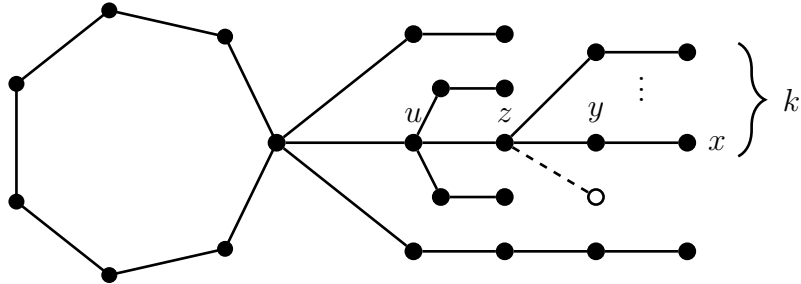
\begin{figure}[H]
	\begin{tikzpicture}
		[line width = 1pt, scale=1.2,
		empty/.style = {circle, draw, fill = white, inner sep=0mm, minimum size=2mm}, full/.style = {circle, draw, fill = black, inner sep=0mm, minimum size=2mm}, full1/.style = {circle, draw, fill = black, inner sep=0mm, minimum size=2mm}]
		
		\foreach \x in {1,...,7} {
			\draw ({360/7 * (\x - 1)}:1.5) -- ({360/7 * \x}:1.5);
			\fill ({360/7 * \x}:1.5) circle (2.5pt);
		}
		
		\node [full1](a1) at (3,1.2) {};
		\node [full1](a2) at (3,-1.2) {};
		\node [full1,label=above:$u$] (u) at (3,0) {};
		\node [full1] (a) at ({360/7 * (7)}:1.5) {};
		\draw (a)--(u);
		\draw (a)--(a1);
		\draw (a)--(a2);
		\node [full1,label=above:$z$] (z) at (4,0) {};
		\node [full1](a12) at (4,1.2) {};
		\node [full1](a22) at (4,-1.2) {};
		\node [full1](a23) at (5,-1.2) {};
		\node [full1](a24) at (6,-1.2) {};
		\node [full1](u1) at (3.3,0.6) {};
		\node [full1](u2) at (3.3,-0.6) {};
		\node [full1](u12) at (4,0.6) {};
		\node [full1](u22) at (4,-0.6) {};
		\draw (u)--(z);
		\draw (a1)--(a12);
		\draw (a2)--(a22);
		\draw (a23)--(a22);
		\draw (a23)--(a24);
		\draw (u)--(u1);
		\draw (u)--(u2);
		\draw (u1)--(u12);
		\draw (u2)--(u22);
		
		\node [label=above:$\vdots$] (w) at (5.5,0.2) {};
		\node [full1,label=above:$y$] (y) at (5,0) {};
		\node [full1,label=right:$x$] (x) at (6,0) {};
		\node [full1](z1) at (5,1) {};
		\node [full1](z12) at (6,1) {};
		\node [empty](z2) at (5,-0.6) {};
		\draw (z)--(z1);
		\draw[dashed] (z)--(z2);
		\draw (z1)--(z12);
		\draw (z)--(y);
		\draw (x)--(y);
		
			\draw[decorate,decoration={brace,amplitude=10pt},xshift=-4pt,yshift=0pt]
		(6.7,1.1) -- (6.7,-0.15) node [black,midway,xshift=0.7cm] {$k$};	
		
	\end{tikzpicture}
	\caption{Subcase 2.3. $d\geq 3$}
	\label{fig7}
\end{figure}

In this subcase, $d_G(y)=2$, and another adjacent vertex of $y$, denoted as $z$, is not in $C$. It is easy to see that the vertex $z$, apart from one adjacent vertex, say $u$, is connected to some pendant edges and at most one leaf. Let's assume $z$ is connected to $k$ pendant edges and $s$ leaves. Clearly $k\geq 1$ and $s\leq 1$. See Figure \ref{fig7}.

We have
\[mis(G)=mis(G-\{x,y\})+mis(G-\{x,y,z\}).\]
Moreover, $G-\{x,y\}$ is a unicyclic graph of order $n-2$ and independence number $\lfloor \frac{n}{2}\rfloor-1$.

If $s=1$, then $G-\{x,y,z\}$ is a disjoint union of a unicyclic graph $G_1$ and $k-1$ isolated edges, along with one isolated vertex. Furthermore, $G_1$ is a unicyclic graph of order $n-(2k+2)$ and independence number $\lfloor\frac{n}{2}\rfloor-(k+1)$. 
If $n-(2k+2)=3$, then by direct calculation,
\[mis(G)=3\cdot 2^k+2\geq \ell(k+3)=\ell(\frac{n+1}{2}).\]
If $n-(2k+2)\geq 5$, then by inductive hypothesis,
\begin{align*}
mis(G)&=mis(G-\{x,y\})+mis(G-\{x,y,z\})\\
&\geq \ell(\frac{n-1}{2})+\ell(\frac{n-(2k+2)+1}{2})\cdot 2^{k-1}\\
&\geq \ell(\frac{n-1}{2})+\ell(\frac{n-3}{2}) \ \text{(by Lemma \ref{lem10})}\\
&=\ell(\frac{n+1}{2}).
\end{align*}

If $s=0$, then $G-\{x,y,z\}$ is a disjoint union of a unicyclic graph $G_1$ and $k-1$ isolated edges. Furthermore, $G_1$ is a unicyclic graph of order $n-(2k+1)$ and independence number at most $\lfloor\frac{n}{2}\rfloor-k$. 

Let $n'=n-(2k+1)$. Then $\alpha(G_1)\leq \lfloor\frac{n}{2}\rfloor-k=\frac{n'}{2}$. On the other hand, since $G_1$ is a unicycle of order $n'$, 
\[\alpha(G_1)\geq\frac{n'}{2},\]
which implies that $\alpha(G_1)=\frac{n'}{2}$.

We show that $u\notin V(C)$. Suppose, to the contrary, that $u\in V(C)$. Then,  
\[\alpha(G_1-u)\geq\lceil\frac{n'-1}{2}\rceil=\frac{n'}{2} \ (\text{since $G_1-u$ is a forest of order $n'-1$}),\]
and
\[\alpha(G)\geq \alpha(G_1-u)+(k+1)\geq \frac{n'}{2}+(k+1)=\frac{n+1}{2}>\lfloor\frac{n}{2}\rfloor,\]
which is a contradiction.

Now, $G_1-u$ is a disjoint union of a unicyclic graph $G'_1$ and a forest $F$. Let $n_1=|V(G'_1)|$, $\alpha_1=\alpha(G'_1)$, $n_2=|V(F)|$, and $\alpha_2=\alpha(F)$. Clearly, $n_1+n_2=n'-1$ ($n'-1$ is odd). 

Since
\begin{align*}
\alpha(G)&=\lfloor\frac{n}{2}\rfloor=\frac{n_1+n_2+2k+1}{2}=\frac{n_1+n_2-1}{2}+(k+1)\\
&\geq \alpha(G_1-u)+(k+1)\\
&=\alpha_1+\alpha_2+(k+1)\\
&\geq \lfloor\frac{n_1}{2}\rfloor+\lceil\frac{n_2}{2}\rceil+(k+1),
\end{align*}
we have that $n_1$ is odd and $n_2$ is even, and $\alpha_1=\lfloor\frac{n_1}{2}\rfloor$, and $\alpha_2=\frac{n_2}{2}$.

If $n_1=3$ and $n_2=0$, then by direct calculation,
\[mis(G)=3\cdot 2^k+2\geq \ell(k+3)=\ell(\frac{n+1}{2}).\]
If $n_1=3$ and $n_2\geq2$, then
\begin{align*}
	mis(G)&=mis(G-\{x,y\})+mis(G-\{x,y,z\})\\
	&\geq mis(G-\{x,y\})+mis(G-\{x,y,z,u\})\\
	&\geq \ell(\frac{n-1}{2})+3\cdot g(n_2-\frac{n_2}{2})\cdot 2^{k-1}\\
	&\geq\ell(\frac{n-1}{2})+\ell(\frac{n_2}{2}+2)\cdot 2^{k-1}\\
	&=\ell(\frac{n-1}{2})+\ell(\frac{n-3}{2}) \ \text{(by Lemma \ref{lem10})} \\
	&=\ell(\frac{n+1}{2}).
\end{align*}
If $n_1\geq 5$, then by inductive hypothesis,
\begin{align*}
	mis(G)&=mis(G-\{x,y\})+mis(G-\{x,y,z\})\\
	&\geq mis(G-\{x,y\})+mis(G-\{x,y,z,u\})\\
	&\geq \ell(\frac{n-1}{2})+\ell(\frac{n_1+1}{2})\cdot g(n_2-\frac{n_2}{2})\cdot 2^{k-1}\\
	&\geq \ell(\frac{n-1}{2})+\ell(\frac{n-3}{2})\ \text{(by Lemmas \ref{lem6} and \ref{lem10})}\\
	&=\ell(\frac{n+1}{2}).
\end{align*}

Now, we have proved that Inequality (\ref{inequ1}) holds.\\

Next, we show that Inequality (\ref{inequ1}) is sharp.

When $n\leq 5$, it is easy to check that the inequality is sharp. We suppose that $n\geq 6$. 

When $n-\alpha=2$, a unicyclic graph $G$ of order $n$, which is constructed by connecting one vertex of a triangle to $n-3$ isolated vertices, satisfies that $\alpha(G)=n-2$ and $mis(G)=3=t(n,\alpha)$. 

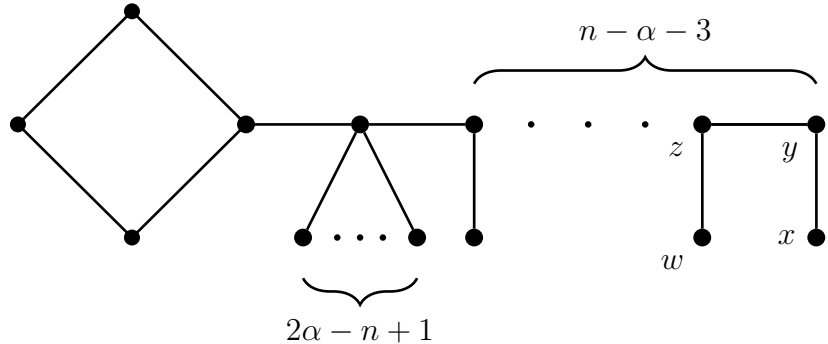
\begin{figure}[H]
		\begin{tikzpicture}
			[line width = 1pt, scale=1.5,
			empty/.style = {circle, draw, fill = white, inner sep=0mm, minimum size=2mm}, full/.style = {circle, draw, fill = black, inner sep=0mm, minimum size=2mm}, full1/.style = {circle, draw, fill = black, inner sep=0mm, minimum size=0.5mm}]
			\foreach \x in {1,2,3,4} {
				\draw ({360/4 * (\x - 1)}:1) -- ({360/4 * \x}:1);
				\fill ({360/4 * \x}:1) circle (2pt);
			}
			\node [full](x) at ({360/4 * 4}:1) {};
			
			\node [full](b) at (2,0) {};
			\node [full](b1) at (1.5,-1) {};
			\node [full1](b2) at (1.8,-1) {};
			\node [full1](b3) at (2,-1) {};
			\node [full1](b4) at (2.2,-1) {};
			\node [full](b5) at (2.5,-1) {};
			\node [full](c) at (3,0) {};
			\node [full](c1) at (3,-1) {};
			\node [full, label=below left:$z$](d) at (5,0) {};
			\node [full, label=below left:$w$](d1) at (5,-1) {};
			\node [full, label=below left:$y$](e) at (6,0) {};
			\node [full, label=left:$x$](e1) at (6,-1) {};
		
			\node [full1](u) at (3.5,0) {};
			\node [full1](v) at (4,0) {};
			\node [full1](w) at (4.5,0) {};
			
			\draw (b) -- (x);
			\draw (b) -- (c);
			\draw (d) -- (e);
			\draw (b) -- (b1);
			\draw (b) -- (b5);
			\draw (c) -- (c1);
			\draw (d) -- (d1);
			\draw (d) -- (e);
			\draw (e) -- (e1);
		
			\draw[decorate,decoration={brace,amplitude=10pt},xshift=0pt,yshift=-4pt]
			(3,0.5) -- (6,0.5) node [black,midway,yshift=0.7cm] {$n-\alpha-3$};
			\draw[decorate,decoration={brace,amplitude=10pt},xshift=0pt,yshift=4pt]
			(2.5,-1.5) -- (1.5,-1.5) node [black,midway,yshift=-0.7cm] {$2\alpha-n+1$};
		\end{tikzpicture}
		\caption{$H_{n,\alpha}$}
		\label{fig8}
	\end{figure}

Consider the case when $\lceil \frac{n}{2} \rceil \leq \alpha < n-2$. A unicyclic graph of order $n$ and independence number $\alpha$ that is shown in Figure \ref{fig8} is denoted by $H_{n,\alpha}$. We prove by induction on $n-\alpha$ that $mis(H_{n,\alpha})=h(n-\alpha)$.

A few simple calculations show that $mis(H_{n,\alpha})=4=h(3)$ when $n-\alpha=3$, and $mis(H_{n,\alpha})=6=h(4)$ when $n-\alpha=4$.

Assume that $n-\alpha\geq 5$, and for integer $3\leq k\leq n-1$, if $n-\alpha=k$, $mis(H_{n,\alpha})=h(k)$. Now, by inductive hypothesis, 
\begin{align*}
	mis(H_{n,\alpha})&=mis(H_{n,\alpha}-\{x,y\})+mis(H_{n,\alpha}-\{x,y,z,w\})\\
	&=mis(H_{n-2,\alpha-1})+mis(H_{n-4,\alpha-2})\\
	&=h(n-\alpha-1)+h(n-\alpha-2)=h(n-\alpha).
\end{align*}

Consider the case when $n\geq 7$, and $n$ is odd, and $\alpha = \lfloor \frac{n}{2} \rfloor$.
A unicyclic graph of order $n$ and independence number $\alpha=\lfloor \frac{n}{2} \rfloor$ that is shown in Figure \ref{fig9} is denoted by $L_n$. 
We prove by induction on $n$ that $mis(L_n)=\ell(\frac{n+1}{2})$.
 
A few simple calculations show that $mis(L_{7})=7=\ell(4)$ and $mis(L_{9})=12=\ell(5)$.
 
Assume that $n\geq 11$, and for any odd integer $7\leq k\leq n-2$, $mis(L_{k})=\ell(\frac{k+1}{2})$. By inductive hypothesis, 
\begin{align*}
	mis(L_n)&=mis(L_n-\{x,y\})+mis(L_n-\{x,y,z,w\})\\
	&=\ell(\frac{n-1}{2})+\ell(\frac{n-3}{2})\\
	&=\ell(\frac{n+1}{2}).
\end{align*}

	\begin{figure}[H]
		\begin{tikzpicture}
			[line width = 1pt, scale=1.5,
			empty/.style = {circle, draw, fill = white, inner sep=0mm, minimum size=2mm}, full/.style = {circle, draw, fill = black, inner sep=0mm, minimum size=2mm}, full1/.style = {circle, draw, fill = black, inner sep=0mm, minimum size=0.5mm}]
			\foreach \x in {1,...,7} {
				\draw ({360/7 * (\x - 1)}:1) -- ({360/7 * \x}:1);
				\fill ({360/7 * \x}:1) circle (2pt);
			}
			\node [full](x) at ({360/7 * 7}:1) {};
			
			\node [full](b) at (2,0) {};
			\node [full](b1) at (2,-1) {};
			\node [full](c) at (3,0) {};
			\node [full](c1) at (3,-1) {};
			\node [full,label=, label=below left:$z$](d) at (5,0) {};
			\node [full,label=, label=below left:$w$](d1) at (5,-1) {};
			\node [full, label=, label=below left:$y$](e) at (6,0) {};
			
			\node [full, label=left:$x$](e6) at (6,-1) {};
			\node [full1](u) at (3.5,0) {};
			\node [full1](v) at (4,0) {};
			\node [full1](w) at (4.5,0) {};
			
			\draw (b) -- (x);
			\draw (b) -- (c);
			\draw (d) -- (e);
			\draw (b) -- (b1);
			\draw (c) -- (c1);
			\draw (d) -- (d1);
			\draw (e) -- (e6);
			\draw[decorate,decoration={brace,amplitude=10pt},xshift=0pt,yshift=-4pt]
			(2,0.5) -- (6,0.5) node [black,midway,yshift=0.7cm] {$\frac{n-7}{2}$};
		\end{tikzpicture}
		\caption{$L_n$}
		\label{fig9}
	\end{figure}

We complete the proof of Theorem \ref{thm2}. \qed
	

\section*{Use of AI tools declaration}
The authors declare that they have not used Artificial Intelligence (AI) tools in the creation of this article.

\section*{Data Availability}
Data sharing is not applicable to this paper as no datasets were generated or analyzed during
the current study.

\section*{Acknowledgments}
This work was supported by Beijing Natural Science Foundation (No. 1232005). 

\section*{Conflict of interest}
The authors declare that they have no conflicts of interest.

	\bibliographystyle{unsrt}
	
\end{document}